\documentclass[12pt]{article}
\usepackage{amsfonts}
\usepackage{mathrsfs}
\usepackage{amsmath}
\usepackage{latexsym}
\usepackage{amssymb}
\usepackage{amsthm}
\usepackage{amscd}

\title{Universal Enveloping Algebras of
Braided m-Lie Algebras }
\author{  \small
Lingwei Guo,  Shouchuan Zhang ,  Jieqiong He  \\
\small $a$. Department  of Mathematics,
Hunan University\\  \small   Changsha  410082, \ P.R. China \\
\small $b$. School of Mathematics and Physics, The University of Queensland\\
\small Brisbane 4072, Australia\\
 }

\begin{document}
\newtheorem{Proposition}{\quad Proposition}[section]
\newtheorem{Theorem}[Proposition]{\quad Theorem}
\newtheorem{Definition}[Proposition]{\quad Definition}
\newtheorem{Corollary}[Proposition]{\quad Corollary}
\newtheorem{Lemma}[Proposition]{\quad Lemma}
\newtheorem{Example}[Proposition]{\quad Example}
\newtheorem{Remark}[Proposition]{\quad Remark}

\maketitle \addtocounter{section}{-1}

\numberwithin{equation}{section}

\date{}

\begin {abstract} Universal enveloping algebras of braided
m-Lie algebras and PBW theorem are obtained by means of
Combinatorics on words.


keywords: Braided Lie algebra, Universal enveloping algebras.
\end {abstract}

\section{Introduction}\label {s0}

The theory of Lie superalgebras has been developed systematically,
which includes  the representation theory and classifications of
simple Lie superalgebras and their varieties \cite {Ka77} \cite
{BMZP92}.  In many physical applications or in pure mathematical
interest, one has to consider not only ${\bf Z}_2$- or ${\bf Z}$-
grading but also $G$-grading of Lie algebras, where $G$ is an
abelian group equipped with a skew symmetric bilinear form given by
a 2-cocycle. Lie algebras in symmetric and more general categories
were discussed in \cite {Gu86} and \cite {GRR95}. A sophisticated
multilinear version of the Lie bracket was considered in \cite
{Kh99} \cite {Pa98}. Various generalized Lie algebras have already
appeared  under different names, e.g. Lie color algebras, $\epsilon
$ Lie algebras \cite {Sc79}, quantum and braided Lie algebras,
generalized Lie algebras \cite {BFM96} and $H$-Lie algebras \cite
{BFM01}.

In  \cite {Ma94c}, Majid introduced braided Lie algebras from
geometrical point of view, which have attracted attention in
mathematics and mathematical physics (see e.g. \cite {Ma95b} and
references therein).

In paper \cite {ZZ04}, braided m-Lie algebras was introduced,  which
generalize Lie algebras, Lie color algebras and quantum Lie
algebras. Two classes of braided m-Lie algebras are given, which are
generalized matrix braided m-Lie algebras and braided m-Lie
subalgebras of $End _F M$, where $M$ is a Yetter-Drinfeld module
over $B$ with dim $B< \infty $ . In particular, generalized
classical braided m-Lie algebras $sl_{q, f}( GM_G(A),  F)$  and
$osp_{q, t} (GM_G(A), M, F)$  of generalized matrix algebra
$GM_G(A)$ are constructed  and their connection with special
generalized matrix Lie superalgebra $sl_{s, f}( GM_{{\bf Z}_2}(A^s),
F)$  and orthosymplectic generalized matrix Lie super algebra
$osp_{s, t} (GM_{{\bf Z}_2}(A^s), M^s, F)$  are established. The
relationship between representations of braided m-Lie algebras and
their  associated algebras are established. In \cite{WZZ15}, the relationship between Nichols braided Lie algebras
 and  Nichols algebras  are studied and provided a new method to determine when a Nichols algebra is finite dimensional.
It was shown that Nichols algebra $\mathfrak B(V)$ is
finite-dimensional if and only if Nichols braided Lie algebra $\mathfrak L(V)$ is finite-dimensional
if there does not exist any $m$-infinity element in $\mathfrak B(V)$.
In \cite{WZZ16} we will strengthen this result by showing that the condition
``there does not exist any $m$-infinity element in $\mathfrak B(V)$" can be dropped.

In this paper we follow
paper \cite {ZZ04} and obtain  universal enveloping algebras of
braided m-Lie algebras and PBW theorem by means of Combinatorics on
words.

Throughout, $ F$ is a field,

\section {Braided  m-Lie Algebras}

 We recalled two concepts.

\begin {Definition} \label {4'.1.1} (See \cite {ZZ04})
Let  $(L, [\ \ ])$ be an object in the braided tensor category $
({\cal C }, C)$ with morphism $[\ \ ] : L \otimes L \rightarrow L$.
If there exists an algebra $(A, m)$ in  $ ({\cal C }, C)$ and
monomorphism $\phi : L \rightarrow A$ such that   $\phi [\  \ ] = m
(\phi \otimes \phi ) - m (\phi \otimes \phi )  C_{L, L},$   then
$(L, [ \ \  ])$ is called a braided m-Lie algebra in $ ({\cal C },
C)$ induced by multiplication of $A$ through $\phi $. Algebra $(A,
m) $ is called an algebra  associated to $(L, [ \ \ ])$.

\end {Definition}
A Lie  algebra is a braided m-Lie algebra in the category of
ordinary vector spaces, a Lie color algebra is a braided  m-Lie
algebra
 in symmetric braided tensor category $ ({\cal M} ^{FG}, C^r)$ since the
 canonical map $\sigma: L \rightarrow U(L)$ is injective (see
\cite [Proposition 4.1]{Sc79}), a  quantum Lie  algebra is a braided
m-Lie algebra in the Yetter-Drinfeld category $ (^B_B{\cal YD}, C)$
by \cite [Definition 2.1 and Lemma 2.2]{GM03}), and a ``good"
braided Lie  algebra is a braided  m-Lie algebra
 in the Yetter-Drinfeld category $ (^B_B{\cal YD}, C)$ by
 \cite [Definition 3.6 and Lemma 3.7]{GM03}).
For  a cotriangular Hopf algebra $(H, r)$, the $(H,r)$-Lie algebra
defined in \cite [4.1] {BFM01} is a
 braided m-Lie  algebra in the braided
tensor category $({}^H{\cal M}, C^r)$.  Therefore, the braided m-Lie
algebras generalize most known generalized Lie algebras.

For an algebra $(A, m)$ in $({\cal C}, C)$, obviously $L = A$ is  a
braided m-Lie algebra under operation $ [\ \ ] = m  - m   C_{L, L}$,
which  is induced by $A$ through $id _A$. This  braided  m-Lie
algebra is written as $A^-$.

\begin {Definition}\label {1} (see \cite {Zh99})
Let H be a Hopf algebra, if there exist the left action $(V, \alpha)$ and left coaction $(V,  \delta )$ of $(V, \alpha, \delta)$,
satisfy:$\forall v\in V$
 $$\delta(\alpha(h\otimes v))=\delta(h.v)=\sum h_{(1)}v_{(-1)}sh_{(3)}\otimes h_{(2)}.v_{0}$$
then H is a Yetter-Drinfeld module on $(V, \alpha, \delta )$. all of
the Yetter-Drinfeld $H$-module construct the braided tensor category
, called Yetter-Drinfeld module category,denote as $^H_H {\cal YD}$.
\end {Definition}

\begin {Definition}\label {1}
Suppose G is an group, k a field,if $\chi :G\longrightarrow k-0$ is
a homomorphism on multiplication group , then $\chi$ is character of
$G$.all character of $G$ construct a group, called the character
group of G, denote as $\hat{G}.$
\end {Definition}

If $L$ is the direct sum of one dimensional Yetter-Drinfeld $kG$
module,then $L$ is called pointed Yetter-Drinfeld $kG$ module.
Further,if $L$ is a Yetter-Drinfeld $kG$ Lie algebra£¬ then  $L$ is
called pointed Yetter-Drinfeld $kG$ Lie algebra.

\begin {Lemma}\label {1'} (see \cite [Lemma 2.2]{ZZC04})
$L$ is pointed YD- Lie algebra on $kG$,if and only if there exist a
set of homogeneous basis of $L$: $\{ x_i \mid i\in J\}$ and a set of
character:$\{ \chi_i \mid i\in J\}\subseteq \hat G$ ,for any $i\in
J, h \in G $,such that $h\cdot x_i = \chi_i (h) x_i$
\end {Lemma}

\section{Jacobi identity}

Let $G$ be a group and $\chi $ a bicharacter of $G$, i.e. $\chi$ is
a map from $G \times G$ to $F$ satisfying  $\chi (ab, c) = \chi (a,
c)\chi (b, c)$, $\chi (a, bc) = \chi (a, b)\chi (a, c) $  and $\chi
(a, e) =1 = \chi (e, a)$ for any $a, b, c \in G$, where $e$ is the
unit element of $G$.

If $V$ is a $G$-graded vector space and $h\cdot x = \chi (\mid x
\mid, h) = x$ for any $h\in G$, homogeneous element $x\in V$. It is
clear that $(V, \delta, \alpha)$ is a $kG$-{\rm YD} module. In this
case, we have $C(b\otimes c)=\chi(b,c)c\otimes b$ for any
homogeneous elements  $a, b \in V$.

For an algebra $(A, m)$ in $(^{kG}_{kG}{\cal YD}, C)$,  $L \subseteq
A$ is a braided m-Lie algebra under operation $ [\ \ ] = m  - m
C_{L, L}$, i.e.  $[x  y]=x  y-\chi(x,  y)y  x$ for any homogeneous
elements $x, y \ L.$

An $A$-module $M$ is called {\it pointed} if $M=0$ or $M$ is a
direct sum of one dimensional $A$-modules.

\begin {Lemma} (See \cite {Kh99}) If $L$ is braided
m-Lie algebra, then Jacobi identity holds:
 \begin{eqnarray}\label {Jacobi} [[a  b]  c]-[a  [b
c]]+\chi(a,b)b [a c]-\chi(b,c)[a c] b=0 \end {eqnarray} for any
homogeneous elements $a, b, c \in L$.
\end {Lemma}
{\bf Proof.}
\begin{eqnarray*}
 \mbox {the left side }&=& abc-\chi(a,b)bac-\chi(ab,c)cab+\chi(ab,c)\chi(a,b)cba\\
& &  -abc+\chi(b,c)acb+\chi(a,bc)bca-\chi(a,bc)\chi(b,c)cba \\
& &  +\chi(a,b)bac-\chi(a,b)\chi(a,c)bca-\chi(b,c)acb+\chi(b,c)\chi(a,c)cab\\
&=&  0. \Box
\end{eqnarray*}

\begin {Lemma}
If $(A,\delta,\alpha)$ is a algebra of pointed $^KG_KG{\cal
YD}$,then for any  homogeneous elements $x,y\in V$, there exist a
bicharacter on $G$,such that $C(x,y)=\chi(x,y)y\otimes x$
\end {Lemma}
{\bf Proof.} because A is pointed $^KG_KG{\cal YD}$ module,for any
homogeneous element $x\in A$,there exist a algebra homomorphism $
\chi_{x}: A\rightarrow F$,such that
\begin{eqnarray*}
\chi_{x}(h)x&=&h\cdot x\\
h\cdot (xy)&=&(h\cdot x)(h\cdot y)
\end{eqnarray*}
First we define $\chi: G\times G\rightarrow F$,such that
$\chi(h,x)=\chi_{x}(h)$,then $\chi$ is a bicharacter.

i) $\chi(c,ab)=\chi(c,a)\chi(c,b)$.

Actually for $\forall x\in V_{a},y\in V_{b},z\in V_{c}$,On the one
hand $h\cdot (xy)=(h\cdot x)(h\cdot y)=\chi_{a}(h)\chi_{b}(h)xy$,on
the other hand $h\cdot (xy)=\chi_{ab}(h)xy$,so
$\chi_{a}(h)\chi_{b}(h)=\chi_{ab}(h)$.

ii)$\chi(ab,c)=\chi(a,c)\chi(b,c)$

Actually
$\chi(ab,c)=\chi_{c}(ab)=\chi_{c}(a)\chi_{c}(b)=\chi(a,c)\chi(b,c)$

iii) $\chi(e,a)=1$. because $\chi(e,a)=\chi_{a}(e)$,$\chi_{a}$is a
homomorphism from group G onto F,then $\chi_{a}(e)=1$.

iv) $\chi(a,e)=1$.because A is algebra of a KG-module,then $1_{A}\in
A_{e}$,and $a\cdot 1_{A}=\epsilon(a)1_{A}=1_{A}$,$a\cdot
1_{A}=\chi_{e}(a)1_{A}$,then $\chi_{e}(a)1_{A}=1_{A}$,that is
$\chi_{e}(a)=1$

\section{ Universal enveloping algebras of braided
m-Lie algebras and PBW theorem    }

Let $E$ be a homogeneous  basis of braided  m-Lie algebra $L$ and
$B$ a set. Let $B^*$ denote the set of all words on $B$ and
 $\varphi$  a bijective map from $E$ to $B$. Define  $[bc]=\varphi([ef])$ for any
  $b = \varphi(e)$, $c = \varphi(f)$, $e, f \in E$.

  Let $\prec$ be an order of $B$ and $P = : \{b_{1}b_{2}\cdots  b_{n}\ | \  b_{i}\in B,
 b_n \prec b_{n-1} \prec \cdots \prec b_1, n \in \mathbb N \}$.

 For any $w\in B^*$, let $\nu (w)$ denote the number of  elements in
 set $\{ (r, s, t) \ | \ w=rasbt;  a,b \in B, r, s, t \in B^*, a\prec
 b,\}$ and $\nu(w)$ is called the index of $w$.
 Obviously, we have
\[
v(ubav)=v(uabv)-1
\] for any $a,b \in B, u, v  \in B^*, a\prec
 b$. We also have that $\nu (w)=0$ if and only if $w\in F$.

\begin {Lemma} \label {3.1}
 There exists $\lambda : F[B] \rightarrow FP$ such that

(i) $\lambda(f)= f,$ $f\in P$;

(ii)  $\lambda(ubcv)=\chi(b,c)\lambda(ucbv)+\lambda(u[bc]v),u,v\in
B^{\ast},b,c\in B $;

(iii)  $\lambda(uv)=\lambda(\lambda(u)v)=\lambda(u\lambda(v)),u,v\in
F[B]$.
\end {Lemma}
{\bf Proof.} For $w \in B^*$, we define $\lambda (w)$ using an
induction first on the length and second on the index. If $w \in B$,
 define $\lambda (w) = w$. Let the length of  $w$ be larger than 1 and define
 $ \lambda (w) =: \chi(b,c)\lambda(ucbv)+\lambda(u[bc]v)$ for $w=ubcv$ with
  $b, c\in B$, $u, v \in B^*$. Now we show that the definition is well-defined. For
  $w=ubcv=u'b'c'v'$ with
  $b, c, b', c'\in B$, $u, v, u', v' \in B^*$, we only need  show that
   \begin {eqnarray} \label {e3.1}\chi(b,c)\lambda(ucbv)+\lambda(u[bc]v) =
 \chi(b',c')\lambda(u'c'b'v')+\lambda(u'[b'c']v'). \end {eqnarray}

We show this by following two  steps.

 ($1^\circ$) If  $|u|\preceq |u'|-2$, then $u'=ubct,v=tb'c'v',t\in B^{*}$.
 By induction hypothesis we have
\begin{eqnarray*}
\mbox {the left side  }  &=& \chi(b,c)\chi(b',c')\lambda(ucbtc'b'v')+\chi(b,c)\lambda(ucbt[b'c']v')\\
& &  +\chi(b',c')\lambda(u[bc]tc'b'v')+\lambda(u[bc]t[b'c']v')
\end{eqnarray*}
and
\begin{eqnarray*}
\mbox {the right side  }  &=& \chi(b,c)\chi(b',c')\lambda(ucbtc'b'v')+\chi(b,c)\lambda(ucbt[b'c']v')\\
& &  +\chi(b',c')\lambda(u[bc]tc'b'v')+\lambda(u[bc]t[b'c']v').
\end{eqnarray*}
Thus  (\ref {e3.1}) holds.

($2^\circ$) If $|u|=|u'|-1$£¬then $u'=ub,c=b',v=c'v'$. We only need
show $
\chi(a,b)\lambda(rbacs)+\lambda(r[ab]cs)=\chi(b,c)\lambda(racbs)+\lambda(ra[bc]s)
$. By induction hypothesis we have
\begin{eqnarray*}
\mbox {the left side  } &=&\chi(a,b)\{\chi(a,c)\lambda(rbcas)+\lambda(rb[ac]s)\}+\lambda(r[ab]cs)\\
&=&\chi(a,b)\chi(a,c)\lambda(rbcas)+\chi(a,b)\lambda(rb[ac]s)+\lambda(r[ab]cs)\\
&=&\chi(a,b)\chi(a,c)\{chi(b,c)\lambda(rcbas)+\lambda(r[bc]as)\}+\chi(a,b)\lambda(rb[ac]s)+\lambda(r[ab]cs)\\
&=&\chi(a,b)\chi(a,c)chi(b,c)\lambda(rcbas)+\chi(a,b)\chi(a,c)\lambda(r[bc]as)+\chi(a,b)\lambda(rb[ac]s)\\
& &+\lambda(r[ab]cs)
\end{eqnarray*} and
\begin{eqnarray*}
\mbox {the right side  } &=&\chi(b,c)\{\chi(a,c)\lambda(rcabs)+\lambda(r[ac]bs)\}+\lambda(ra[bc]s)\\
&=&\chi(b,c)\chi(a,c)\lambda(rcabs)+\chi(b,c)\lambda(r[ac]bs)+\lambda(ra[bc]s)\\
&=&\chi(b,c)\chi(a,c)\{\chi(a,b)\lambda(rcbas)+\lambda(rc[ab]s)\}+\chi(b,c)\lambda(r[ac]bs)+\lambda(ra[bc]s)\\
&=&\chi(b,c)\chi(a,c)\chi(a,b)\lambda(rcbas)+\chi(b,c)\chi(a,c)\lambda(rc[ab]s)+\chi(b,c)\lambda(r[ac]bs)\\
& &+\lambda(ra[bc]s).
\end{eqnarray*}
Thus
\begin{eqnarray*}
&&\mbox {the left side  } - \mbox {the right side  }\\
 &=&\{\chi(a,b)\chi(a,c)\lambda(r[bc]as)-\lambda(ra[bc]s)\}+\{\chi(a,b)\lambda(rb[ac]s)-\chi(b,c)\lambda(r[ac]bs)\}\\
& &+\{\lambda(r[ab]cs)-\chi(b,c)\chi(a,c)\lambda(rc[ab]s)\}\\
&=&-\lambda(r[a[bc]]s)+\lambda(\chi(a,b)rb[ac]s-\chi(b,c)r[ac]bs)+\lambda(r[[ab]c]s)\\
&=&0 \ \ \ \ { \mbox {(by Jacobi identity)}}.
\end{eqnarray*}

For (iii), we  use an induction first on the length and second on
the index.  Assume   $|w_1|\neq |w|$ and  $w=w_{1}w_{2}$. If
$w_{1}=ubct$, $b, c\in B$, $u, t \in B^*$, then
\begin{eqnarray*}
&\lambda(w)&=\chi(b,c)\lambda(ucbtw_{2})+\lambda(u[bc]tw_{2})\\
& &=\chi(b,c)\lambda(\lambda(ucbt)w_{2})+\lambda(\lambda(u[bc]t)w_{2})\\
& &=\lambda(\lambda(w_{1})w_{2})
\end{eqnarray*}
If $w_{1}=b$, $w_2 = cv$, $b, c \in B$, $v \in B^*$, then $\lambda
(w) = \lambda ( \lambda (w_1)w_2).$ $\Box$

\begin {Definition}
Suppose $\mathcal{L}$ is a braided  m-Lie algebra in $ ({\cal C },
C)$,$U$ is a algebra,and $\varphi: \mathcal{L}\rightarrow U_{L}$ is
a injective Lie algebra isomorphism,if $T\in  ({\cal C }, C)$,$\psi:
\mathcal{L}\rightarrow T_{L}$ is a Lie algebra homomorphism of $
({\cal C }, C)$,then $\exists \overline{\psi}$
$$ \begin {array} {lcccr} {}&\varphi & {}\\
L& \longrightarrow &U \\
&  \psi  \searrow &  \bar \psi \downarrow  \\
 & & T &  . \end {array}$$
such that $\overline{\psi}\varphi=\psi$.then $U$ is called the
enveloping  braided  algebra of $(\mathcal{L},\varphi)$
\end {Definition}

Let $E$ be a homogeneous  basis of braided  m-Lie algebra $L$ and
$B$ a set. Let $B^*$ denote the set of all words on $B$ and
 $\varphi$  a bijective map from $E$ to $B$. Define  $[bc]=\varphi([ef])$ for any
  $b = \varphi(e)$, $c = \varphi(f)$, $e, f \in E$. Obviously,
  $\varphi$ is a monomorphism from $L$ to $FP,$ the space spanned by
  $P$.

Let $U(L)=: FP$. Define the multiplication of $U$ as follows: $u *
v=\lambda(uv) $ for any $ u, v \in P.$
 By Lemma \ref {3.1} (iii), $U$ is an associative algebra:
$u\ast(v\ast
w)=\lambda(u\lambda(vw))=\lambda(uvw)=\lambda(\lambda(uv)w)=(u\ast
v)\ast w$ for any $u, v, w \in P.$ Obviously, $\lambda $ is an
algebra homomorphism.

\begin {Lemma}
$T(FB)$ is $FG-YD$ module.
\end {Lemma}
{\bf Proof.}
 i) $$ \begin {array} {lcccr} {}&\delta_{T(FB)} & {}\\
T(FB)& \longrightarrow &FB\otimes FB\\
i  \uparrow&  \nearrow (id\otimes i)\delta_{FB}  & \uparrow  \\
 FB&  \longrightarrow &  FG\otimes FB .\end {array}$$

ii) $$ \begin {array} {lcccr} {}&\alpha_{g} & {}\\
FB& \longrightarrow &FB\\
i  \downarrow & \searrow (i\alpha_{g})  & i\downarrow  \\
 T(FB)&  \longrightarrow &  T(B) .\end {array}$$

iii)$\forall g\in G,x_{1},\cdot\cdot\cdot,x_{r}\in B,x_{i}\in
(FB)_{g_{i}}$
\begin {eqnarray*}
\delta(g\cdot(x_{1}\cdot\cdot\cdot
x_{r}))&=&\delta\cdot(\alpha_{g}(x_{1}\cdot\cdot\cdot x_{r}))\\
&=&\delta((g\cdot x_{1})\cdot\cdot\cdot(g\cdot x_{r}))=\delta(g\cdot
x_{1})\cdot\cdot\cdot\delta(g\cdot x_{r})\\
&=&(gg_{1}g^{-1}\otimes (g\cdot
x_{1}))\cdot\cdot\cdot(gg_{r}g^{-1}\otimes (g\cdot x_{r}))\\
&=&(gg_{1}g^{-1})\cdot\cdot\cdot(gg_{r}g^{-1})\otimes
x_{1}\cdot\cdot\cdot x_{r}\\
&=&g(g_{1}\cdot\cdot\cdot g_{r})g^{-1}\otimes x_{1}\cdot\cdot\cdot
x_{r}
\end {eqnarray*}

\begin {Theorem} \label {3.4}(PBW). $(U, \varphi$) is the universal enveloping
algebra of braided Lie algebra $L$.
\end {Theorem}
{\bf Proof.} Define $\bar \psi : FB^* \rightarrow FP$ such that
$\bar \psi \varphi  = \psi$ and $\theta = : \bar \psi \mid _{FP}$,
the restriction of $\bar \psi$ on $FP.$ It is clear that the
following is commutative. $$ \begin {array} {lcccr} {}&\varphi & {}
& \lambda & {}\\
L& \longrightarrow &FB^* &\longrightarrow & FP\\
&  \psi  \searrow &  \bar \psi \downarrow & \swarrow \theta \\
 & & T &  . \end {array}$$
We only need show that $\theta$ is an algebra homomorphism, i.e.
 \begin{eqnarray*}  \theta(r\ast s)&=&
\theta(r) \theta (s) . \end {eqnarray*} We show this by following
several steps.

$(1^\circ)$ If $rs\in P,$ then $\theta(r\ast
s)=\theta(\lambda(rs))=\theta(rs)=\theta(r)\theta(s)$.
  \begin{eqnarray*} \ \ \ \  (2^\circ) \ \ \ \theta(r\ast s)&=&
\theta(\lambda
(rs))\\
&=&\theta ( \lambda (sr \chi (r,s) + [rs]))\\
&=&\theta ( \lambda (sr ))\chi (r,s) + \theta  (\lambda ([rs]))\\
&=&\theta ( s*r )\chi (r,s) + \theta  ([rs]) \ \ \  ( \mbox { the
length of } [rs] <2 \mbox
{ and } \nu (s*r) < \nu (rs) )\\
&=&\theta ( s)\theta (r )\chi (r,s) + \theta  ([rs])\\
&=&\theta ( sr)\chi (r,s) + \theta  ([rs])\\
&=&\theta ( rs) = \theta (r)\theta (s)
\end {eqnarray*}

$(3^\circ)$ If $r=ub,$ $s=cv, u,v\in B^{\ast}$, $b, c\in B, b<c$,
then  \begin{eqnarray*}
\theta(r\ast s)&&=\chi(b,c)\theta(\lambda(ucbv))+\theta(\lambda(u[bc]v))\\
& &=\chi(b,c)\theta((uc) *(bv))+\theta((u[bc]) *v)\\
& &=\chi(b,c)\theta(uc)\theta(bv)+\theta(u[bc])\theta(v)  \ \ \ {
\mbox {(by induction hypothesis)}}\\
& &=\chi(b,c)\theta(u)\theta(cb)\theta(v)+\theta(u)\theta([bc])\theta(v)\\
& &=\theta(u)\theta(bc)\theta(v)\\
& &=\theta(u)\theta(b)\theta(c)\theta(v)\\
& &=\theta(r)\theta(s). \ \ \Box
\end{eqnarray*}


\begin{thebibliography}{BD99}

\bibitem {BFM01} Y. Bahturin, D. Fischman and  S. Montgomery,
Bicharacter, twistings and Scheunert's theorem for Hopf algebra, J.
Alg. {\bf 236} (2001), 246-276.


\bibitem {BFM96}  Y. Bahturin, D. Fischman and S. Montgomery.
On the generalized Lie structure of associative algebras. Israel J.
of Math., {\bf 96}(1996) , 27--48.


\bibitem {BMZP92} Y. Bahturin, D. Mikhalev, M. Zaicev and V. Petrogradsky,
Infinite dimensional Lie superalgebras, Walter de Gruyter Publ.
Berlin, New York, 1992.
\bibitem {GM03} X. Gomez and S. Majid,  Braided Lie algebras and bicovariant
differential calculi over coquasitriangular Hopf algebras, J. Alg.
{\bf 261}(2003), 334--388.

\bibitem {GRR95}  D. Gurevich, A. Radul and V. Rubtsov,
Noncommutative differential geometry related to the Yang-Baxter
equation, Zap. Nauchn. Sem. S.-Peterburg Otdel. Mat. Inst. Steklov.
(POMI) {\bf 199 } (1992); translation in J. Math. Sci. {\bf 77 }
(1995), 3051--3062.

\bibitem {Gu86}  D. I. Gurevich, The Yang-Baxter equation and the
generalization of formal Lie theory, Dokl. Akad. Nauk SSSR, {\bf
288} (1986), 797--801.

\bibitem {Ka77} V. G. Kac. Lie superalgebras. Adv. in Math.,
{\bf 26}(1977) , 8--96.
 \bibitem {Kh99} V. K. Kharchenko,
An existence condition for multilinear quantum operations, J. Alg.
{\bf 217} (1999), 188--228.

\bibitem{Lo83}  M. Lothaire,  Combinatorics on words. London:Cambridge University
Press, 1983.

 \bibitem {Ma93b} S. Majid, Free braided differential calculus, braided binomial
theorem, and the braided exponential map. J. Math. Phys., {\bf 34},
1993, 4843--4856.

\bibitem {Ma94c} S. Majid, Quantum and braided Lie algebras, J. Geom. Phys.
{\bf 13} (1994), 307--356.


  \bibitem {Ma95b} S. Majid, Foundations of  Quantum Group Theory,  Cambradge University Press, 1995.
\bibitem {Sc79}  M. Scheunert. Generalized Lie algebras.
 J. Math. Phys., {\bf 20} (1979), 712--720.

\bibitem {Pa98} B. Pareigis, On Lie algebras in the category of
Yetter-Drinfeld modules. Appl. Categ. Structures,  {\bf 6} (1998),
151--175.

\bibitem{Wo89}  S. L. Woronowicz, Differential calculus on compact matrix
pseudogroups(quantum groups).  Commun. Math. Phys, {\bf 122}(1989)1,
125-170.


\bibitem{ZZ04} S. C. Zhang, Y. Z. Zhang,   Braided m-Lie algebras.
Letters in Mathematical Physics,  2004,  70: 155-167. Also in
math.RA/0308095.


\bibitem[WZZ14] {WZZ14} W. Wu,     S. Zhang and   Y.-Z. Zhang,   Finite dimensional Nichols algebras over finite cyclic groups, J. Lie Theory {\bf 24} (2014),    351-372.

\bibitem[WZZ15] {WZZ15} W. Wu,  S. Zhang and   Y.-Z. Zhang,   Relationship between Nichols braided Lie
algebras and Nichols algebras, J. Lie Theory {\bf 25} (2015),    45-63.

\bibitem[WZZ16] {WZZ16} W. Wu,  S. Zhang and   Y.-Z. Zhang,   On  Nichols braided Lie
algebras. preprint.


\end{thebibliography}
\end {document}